\documentclass[12pt]{amsart}
\newcommand{\sbar}{\ol{\s}}
\makeatletter
\@namedef{subjclassname@2020}{\textup{2020} Mathematics Subject Classification}
\makeatother

\newcommand{\dsum}{\di\sum}

\newcommand{\el}{\ensuremath{\ell}}

\DeclareFontFamily{U}{mathx}{\hyphenchar\font45}
\DeclareFontShape{U}{mathx}{m}{n}{
      <5> <6> <7> <8> <9> <10>
      <10.95> <12> <14.4> <17.28> <20.74> <24.88>
      mathx10
      }{}
\DeclareSymbolFont{mathx}{U}{mathx}{m}{n}
\DeclareFontSubstitution{U}{mathx}{m}{n}
\DeclareMathAccent{\widecheck}{0}{mathx}{"71}

\renewcommand{\kill}[1]{}
\newcommand{\dummy}[1]{\mbox{}}

\makeatletter
\newcommand{\xequal}[2][]{\ext@arrow 0055{\equalfill@}{#1}{#2}}
\def\equalfill@{\arrowfill@\Relbar\Relbar\Relbar}
\makeatother

\renewcommand{\k}{\ensuremath{\ol{\mathrm{P}}}}

\newcommand{\m}{\ensuremath{\infty}}

\renewcommand{\k}[1]{\ensuremath{\left({#1}\right)}}

\newcommand{\bca}{\begin{cases}}
\newcommand{\eca}{\end{cases}}

\newcommand{\mug}{\ensuremath{\infty}}

\newcommand{\dprod}{\di\prod}

\newcommand{\ff}[2]{\ensuremath{\di\fr{#1}{#2}}}

\newcommand{\s}[1]{\ensuremath{\di\int{#1}\,dx}}

\newcommand{\bpic}{\begin{picture}}\newcommand{\epic}{\end{picture}}

\newcommand{\beda}{\begin{edaenumerate}}
\newcommand{\eeda}{\end{edaenumerate}}

\newcommand{\cd}{\cdots}

\newcommand{\q}{\quad}

\newcommand{\bq}{\begin{quote}}\newcommand{\eq}{\end{quote}}

\newcommand{\be}{\begin{enumerate}}\newcommand{\ee}{\end{enumerate}}
\newcommand{\bce}{\begin{center}}\newcommand{\ece}{\end{center}}
\newcommand{\bde}{\begin{description}}\newcommand{\ede}{\end{description}}
\newcommand{\bri}{\begin{flushright}}\newcommand{\eri}{\end{flushright}}
\newcommand{\bb}{\begin{block}}\newcommand{\eb}{\end{block}}
\newcommand{\bt}{\begin{thm}}\newcommand{\et}{\end{thm}}
\newcommand{\bpf}{\begin{proof}}\newcommand{\epf}{\end{proof}}
\newcommand{\bex}{\begin{ex}}\newcommand{\eex}{\end{ex}}
\newcommand{\bexr}{\begin{exr}}\newcommand{\eexr}{\end{exr}}
\newcommand{\bft}{\begin{fact}}\newcommand{\eft}{\end{fact}}
\newcommand{\brk}{\begin{rmk}}\newcommand{\erk}{\end{rmk}}
\newcommand{\ba}{\begin{align*}}\newcommand{\ea}{\end{align*}}
\newcommand{\bexe}{\begin{exe}}\newcommand{\eexe}{\end{exe}}
\newcommand{\tn}{\textnormal}

\newcommand{\bit}{\begin{itemize}}\newcommand{\eit}{\end{itemize}}

\newcommand{\bcm}{}
\newcommand{\ol}{\overline}

\newcommand{\fr}{\frac}

\newcommand{\zz}{\ensuremath{\mathbf{Z}}}

\newcommand{\bd}{\begin{defn}}\newcommand{\ed}{\end{defn}}
\newcommand{\bp}{\begin{prop}}\newcommand{\ep}{\end{prop}}

\newcommand{\eh}{\emph}

\newcommand{\te}{\text}
\newcommand{\wt}{\widetilde}

\newcommand{\di}{\displaystyle}

\renewcommand{\s}{\sigma}

\usepackage[dvipdfmx]{graphicx}
\usepackage[dvipsnames]{xcolor}
\usepackage{float,afterpage}
\usepackage{asymptote,layout}
\usepackage{wrapfig,epic}
\usepackage{amsmath,amssymb,cases,pict2e}

\usepackage{multicol}
\usepackage{tikz}

\usepackage{geometry}
\usepackage{exscale,latexsym,bm}
\usepackage{amssymb,enumerate,amsmath,amsthm,amsfonts}
\usepackage{verbatim,fancybox,wasysym,fancyhdr,type1cm}
\usepackage[frame,all,poly,curve,knot,arrow]{xy}
\usepackage{colortbl}
\usepackage{boxedminipage}
\usepackage{multirow}

\graphicspath{{./figs/}}
\theoremstyle{definition}
\newtheorem{thm}{Theorem}[section]

\newtheorem{lem}[thm]{Lemma}
\newtheorem{prop}[thm]{Proposition}\newtheorem{cor}[thm]{Corollary}

\newtheorem{exr}[thm]{Exercise}

\newtheorem{ex}[thm]{Example}

\newtheorem{defn}[thm]{Definition}\newtheorem{rmk}[thm]{Remark}
\newtheorem{fact}[thm]{Fact}
\newtheorem{block}[thm]{}
\newtheorem*{exe}{Exercise}

\renewcommand{\arraystretch}{1.5}

\newcommand{\tf}{\tfrac}

\begin{document}

\title{New recurrences for divisor sum functions and triangular numbers}
\author{Masato Kobayashi}
\date{\today}                                       

\subjclass[2010]{Primary:11A25;\,Secondary:05A15;}
\keywords{divisor sum functions, 
Jacobi triple product identity, Lambert series, pentagonal numbers, triangular numbers. 
}
\address{Masato Kobayashi\\
Department of Engineering\\
Kanagawa University, 3-27-1 Rokkaku-bashi, Yokohama 221-8686, Japan.}
\email{masato210@gmail.com}


\maketitle
\begin{abstract}
Euler discovered recurrence for 
divisor sum functions as a consequence of the pentagonal numbers theorem. With similar idea 
and also motivated by Ewell's work in 1977, we prove  new recurrences for certain divisor sum functions and  triangular numbers. 
\end{abstract}
\tableofcontents

\section{Introduction: Pentagonal number theorem}

The Euler's \eh{pentagonal number theorem} \cite{euler} claims that, as formal power series, 
\[
\dprod_{n=1}^{\mug} (1-q^{n})
=1-q^{1}-q^{2}+q^{5}+q^{7}-q^{12}-q^{15}+q^{22}+q^{26}-\cd.
\]
Exponents appearing on the right hand side are exactly  \eh{pentagonal numbers}, nonnegative integers in the form
$e_{m}=\tfrac{m(3m-1)}{2}$, $m\in\zz$ (Table \ref{t1}).
{\renewcommand{\arraystretch}{1.25}
\begin{table}[h!]
\caption{pentagonal, triangular, hexagonal numbers}
\label{t1}
\begin{center}
\begin{tabular}{c|ccccccccccccccccccc}
$k$&0&1&2&3&4&5&$\cd$\\\hline
$e_{k}$&0&1&5&12&22&35&$\cd$\\
$e_{-k}$&0&2&7&15&26&40&$\cd$\\\hline
$t_{k}$&0&1&3&6&10&15&$\cd$\\\hline
$h_{k}$&0&1&6&15&28&45&$\cd$\\
$h_{-k}$&0&3&10&21&36&55&$\cd$\\
\end{tabular}
\end{center}
\end{table}%
}

He further proved recurrence for 
the number of partitions $p$ for a natural number $n$:
\[
p(n)=p(n-1)+p(n-2)-p(n-5)-p(n-7)+\cd.
\]
This sum indeed involves only finitely many terms 
with the convention that $p(0)=0$ and $p(m)=0$ whenever $m<0$. 
In history, the great mathematicians, such as Euler, Franklin, Gauss, Hardy, Jacobi, MacMahon, Ramanujan, Rogers, Sylvester, had studied this function. See the book \cite{andrews} by  Andrews and K. Eriksson for details.

Yet there is another recurrence (also by Euler) on the \eh{divisor sum function} (Table \ref{t2}). Let
\[
\s(n)=
\sum_{d|n}d.
\]
Then we have 
\[
\s(n)=\s(n-1)+\s(n-2)-\s(n-5)-\s(n-7)+\cd
\]
where $``\s(0)"$ reads as $n$ (if it appears) and $\s(m)=0$ for $m<0$.
For example, 
\[
\s(7)=\s(6)+\s(5)-\s(2)-``\s(0)"=
12+6-3-7=8,
\]
\[
\s(8)=\s(7)+\s(6)-\s(3)-\s(1)=
8+12-4-1=15.
\]
See 
Thomas-Abdul-Tirupathi \cite{tat} on surprising connections between $p$ and $\s$.

The aim of this article is to prove 
similar results for certain divisor sum functions and  \eh{triangular numbers} as Theorem \ref{th0}  (= Theorem \ref{th1}). 
This is in fact subsequent work of Ewell \cite{ewell}.

\section{Main result}

{\renewcommand{\arraystretch}{1.25}
\begin{table}[h!]
\caption{divisor sum functions}
\label{t2}
\begin{center}
\begin{tabular}{c|ccccccccccccccccccc}
$n$&1&2&3&4&5&6&7&8&9&10&$\cd$\\\hline
$\s(n)$&1&3&4&7&6&12&8&15&13&18&$\cd$\\
$\s_{0}(n)$&0&2&0&6&0&8&0&14&0&12&$\cd$\\
$\s_{1}(n)$&1&1&4&1&6&4&8&1&13&6&$\cd$\\
$\wt{\s}(n)$&1&$-1$&4&$-5$&6&$-4$&8&$-13$&13&$-6$&$\cd$\\
$\sbar(n)$&1&1&4&5&6&4&8&13&13&6&$\cd$\\
\end{tabular}
\end{center}
\end{table}%
}
\begin{defn}
For $n\ge 1$, define 
\[
\s_{0}(n)=
\sum_{
\substack{
d|n\\
d\,\,\text{even}}}d, \q 
\s_{1}(n)=
\sum_{
\substack{
d|n\\
d\,\,\text{odd}}}d
\]
and $\wt{\s}(n)=\s_{1}(n)-\s_{0}(n).$ 
For convenience, 
set $\wt{\s}_{0}(n)=\wt{\s}_{1}(n)=0$ if $n\le 0$.
\end{defn}


A \eh{triangular number} is an integer
 in the form $t_{k}=\tfrac{k(k+1)}{2}$, $k\ge0$.
Let $T=\{t_{k}\mid k\ge0\}$ and \[
\chi_{T}(n)=
\begin{cases}
	1&\text{if }n\in T,	\\
	0&\text{if }n\not\in T.\\
\end{cases}
\]
\begin{thm}
\label{th0}
For all $n\ge 1$, we have 
\[
\wt{\s}(n)=-
\dsum_{k\ge 1}\wt{\s}
\k{n-\ff{k(k+1)}{2}}
+
n\chi_{T}(n).
\]
\end{thm}
For example, 
\[
\wt{\s}(9)=-(\wt{\s}(8)+\wt{\s}(6)+\wt{\s}(3))
+0
=
-(-13-4+4)=
13,
\]
\[
\wt{\s}(10)=-(\wt{\s}(9)+\wt{\s}(7)+\wt{\s}(4))+
10
=-(13+8-5)+10=
-6.
\]
After Lemma \ref{lem1} on Lambert series, we will give a proof of this theorem and discuss its consequence.

\section{Lemma}

In the sequel, we regard all infinite sums and products in $q$ 
as formal power series. 
Notice that we can split 
the infinite product in pentagonal number theorem 
into two parts as 
\[
\dprod_{i=1}^{\mug}(1-q^{i})=
\dprod_{i=1}^{\mug}(1-q^{2i})
(1-q^{2i-1}).
\]
Let us now consider another product
\[
T(q)=\dprod_{i=1}^{\mug} \ff{1-q^{2i}}{1-q^{2i-1}}
\]
with 
\[
T_{0}(q)=
\dprod_{i=1}^{\mug}(1-q^{2i}),\q 
T_{1}(q)=
\dprod_{i=1}^{\mug}\ff{1}{1-q^{2i-1}}
\]
so that $T(q)=T_{0}(q)T_{1}(q)$.
It follows from \eh{Jacobi triple product identity} 
\cite{andrews} that 
\[
T(q)=\dsum_{n=0}^{\mug}q^{n(n+1)/2}.
\]
\begin{lem}\label{lem1}
\[
q\ff{T_{0}'(q)}{T_{0}(q)}=-\dsum_{n=1}^{\mug}\s_{0}(n)q^{n}
,\q 
q\ff{T_{1}'(q)}{T_{1}(q)}=
\dsum_{n=1}^{\mug}\s_{1}(n)q^{n}.
\]
\end{lem}
\begin{proof}
Starting with 
$T_{0}(q)=
\prod_{i=1}^{\mug}(1-q^{2i})$, take its logarithmic derivative and multiply the result by $q$ to obtain 
\[
q\ff{T_{0}'(q)}{T_{0}(q)}
=
- \dsum_{i=1}^{\mug}\ff{2i q^{2i}}{1-q^{2i}}
=
- \dsum_{i=1}^{\mug}
\dsum_{j=0}^{\mug}2i
q^{2i(j+1)}
\]
\[=
- \dsum_{i=1}^{\mug}
\dsum_{k=1}^{\mug}2i
q^{2ik}
=-
\dsum_{n=1}^{\mug}
\sum_{
\substack{
d|n\\
d\,\,\text{even}}}
d
q^{n}
=
-\dsum_{n=1}^{\mug}\s_{0}(n)q^{n}.
\]
Quite similarly, we can show that 
\[
q\ff{T_{1}'(q)}{T_{1}(q)}
=\dsum_{n=1}^{\mug}\s_{1}(n)q^{n}.
\]
\end{proof}

\section{Proof of Theorem}

\begin{thm}\label{th1}
For all $n\ge 1$, we have 
\[
\wt{\s}(n)=-
\dsum_{k\ge 1}\wt{\s}
\k{n-\ff{k(k+1)}{2}}
+n\chi_{T}(n).
\]
\end{thm}
\begin{proof}
Let 
$T_{0}(q), T_{1}(q), T(q)$ as above. 
Thanks to Lemma \ref{lem1}, we have 
\[
q\ff{T'(q)}{T(q)}=
q\ff{T_{0}'(q)}{T_{0}(q)}
+
q\ff{T_{1}'(q)}{T_{1}(q)}=
-
\dsum_{n=1}^{\mug}\s_{0}(n)q^{n}
+
\dsum_{n=1}^{\mug}\s_{1}(n)q^{n}
=\dsum_{n=1}^{\mug}\wt{\s}(n)q^{n},
\]
that is, 
\[
\dsum_{n=0}^{\mug}n \chi_{T}(n)q^{n}
=
\k{\dsum_{\el=0}^{\mug}\chi_{T}(\el)q^{\el}}
\k{\dsum_{m=1}^{\mug}\wt{\s}(m)q^{m}}.
\]
Now equate coefficients of $q^{n}$:
\[
n\chi_{T}(n)=
\dsum_{\el=0}^{n-1}\chi_{T}(\el)\wt{\s}(n-\el).
\]
Rewriting this, we conclude that 
\[
\wt{\s}(n)=-
\dsum_{k\ge 1}\wt{\s}
\k{n-\ff{k(k+1)}{2}}
+n\chi_{T}(n).
\]
\end{proof}
\renewcommand{\ep}{\varepsilon}

\begin{cor}
If $n\not\in T$, then
\[
\sum_{k\ge0}\s_{0}\k{n-\ff{k(k+1)}{2}}
=
\sum_{k\ge0}\s_{1}\k{n-\ff{k(k+1)}{2}}.
\] 
\end{cor}
Observe the case for $n=11$ that 
\[
\s_{0}(11)+\s_{0}(10)+\s_{0}(8)+\s_{0}(5)+\s_{0}(1)
=0+12+14+0+0=26
\]
as well as 
\[
\s_{1}(11)+\s_{1}(10)+\s_{1}(8)+\s_{1}(5)+\s_{1}(1)
=12+6+1+6+1=26.
\]

\section{Consequence}

There is another result (Theorem \ref{th2}) with the idea of 
\eh{hexagonal numbers}.
Since the discussion is analogous,  below we give only its outline. 

Let $h_{m}=m(2m-1)$ $(m\in\zz)$ 
and 
$H=\{h_{m}\mid m\in \zz\}$. 
Notice that $h_{0}=0=t_{0}$ and 
\[
h_{n}=\tf{2n(2n-1)}{2}=t_{2n-1}, \q
h_{-n}=\tf{2n(2n+1)}{2}=t_{2n} 
\]
for $n>0$. 
Thus $H=\{h_{m}\mid m\in\zz\}=\{t_{k}\mid k\ge0\}=T.$ 
It is now natural to introduce the following signed characteristic function:
\[
\ep_{H}(n)=
\begin{cases}
	(-1)^{m}&\te{if } n=h_{m},\\
	0&\tn{ otherwise.}\\
\end{cases}
\]
Further, set 
\[
H(q)=
\dsum_{m=-\mug}^{\mug}
(-1)^{m}q^{h_{m}}=1-q-q^{3}+q^{6}+q^{10}-\cd.
\]
Again, as a consequence of Jacobi triple product identity, one has 
\[
H(q)=
\dprod_{i=1}^{\mug}(1-q^{4i-3})(1-q^{4i-1})(1-q^{4i}).
\]
These factors now suggest us to introduce 
\[
\sbar(n)=
\sum_{
\substack{d|n \\
d\equiv 0, 1, 3 \te{ mod 4}}
}d
\]
and understand $\sbar(n)=0$ whenever $n\le0$.
Then the relation 
\[
q\ff{H'(q)}{H(q)}= \dsum_{n=1}^{\mug}\sbar(n)q^{n}
\]
leads us to the following.
\begin{thm}
\label{th2}
For $n\ge 1$, 
\[
\sbar(n)=
\sum_{
\substack{m\in \zz\\m\ne0}
}
(-1)^{m+1}\sbar\k{
n-m(2m-1)
}
-n\ep_{H}(n).
\]
\end{thm}
For example, the case for $n=10\in H$ is
\[
\sbar(10)=
\sbar(9)+\sbar(7)-\sbar(4)-10
=13+8-5-10=6.
\]
\begin{cor}
If $n\not\in H$, then 
\[
\sum_{j\ge0}\sbar(n-h_{2j})=
\sum_{j\ge0}\sbar(n-h_{2j+1}).
\]
\end{cor}
For example, for $n=11\not\in H$,
\[
\sbar(11)+\sbar(5)+\sbar(1)=
12+6+1=19,
\]
\[
\sbar(10)+\sbar(8)=6+13=19.
\]

%


\begin{center}
\textbf{Acknowledgment.}
\end{center}

\begin{center}
This research arose from Iitaka online seminar in 2020-2022. The author would like to thank the organizer Shigeru Iitaka and other participants. Also Satomi Abe sincerely supported him for his writing.
\end{center}

\end{document}